\documentclass[12pt]{article}

\usepackage{amssymb,fancyheadings}

\newcommand{\nc}{\newcommand}
\nc{\nit}{\noindent}
\nc{\nn}{\nonumber}
\nc{\D}{\partial}
\nc{\diff}[2]{\frac{d #1}{d #2}}
\nc{\diffn}[3]{\frac{d^{#3} #1}{d {#2}^{#3}}}
\nc{\pdiff}[2]{\frac{\partial #1}{\partial #2}}
\nc{\pdiffn}[3]{\frac{\partial^{#3} #1}{\partial{#2}^{#3}}}
\nc{\abs}[1] {\lvert #1 \rvert}
\nc{\cAc}{{\cal A}_c}
\nc{\cE}{{\cal E}}
\nc{\cF}{{\cal F}}
\nc{\cP}{{\cal P}}
\nc{\cV}{{\cal V}}
\nc{\cQ}{{\cal Q}}
\nc{\cGin}{{\cal G}_{\rm in}}
\nc{\cGout}{{\cal G}_{\rm out}}
\nc{\cO}{{\cal O}}
\nc{\Lav}{{\cal L}_{\rm av}}
\nc{\cL}{{\cal L}}
\nc{\cB}{{\cal B}}
\nc{\cZ}{{\cal Z}}
\nc{\cR}{{\cal R}}
\nc{\cT}{{\cal T}}
\nc{\cY}{{\cal Y}}
\nc{\cX}{{\cal X}}
\nc{\cXT}{{{\cal X}(T)}}
\nc{\cBT}{{{\cal B}(T)}}
\nc{\vD}{{\vec \mathcal{D}}}
\nc{\efield}{\mathcal{E}}
\nc{\vE}{{\vec \efield}}
\nc{\vB}{{\vec \mathcal{B}}}
\nc{\vH}{{\vec \mathcal{H}}}
\nc{\ty}{{\tilde y}}
\nc{\tu}{{\tilde u}}
\nc{\tV}{{\tilde V}}
\nc{\Pc}{{\bf P_c}}
\nc{\bx}{{\bf x}}
\nc{\bX}{{\bf X}}
\nc{\bXYZ}{{\bf XYZ}}
\nc{\bY}{{\bf Y}}
\nc{\bF}{{\bf F}}
\nc{\bS}{{\bf S}}
\nc{\dV}{{\delta V}}
\nc{\dE}{{\delta E}}
\nc{\TT}{{\Theta}}
\nc{\dPsi}{{\delta\Psi}}
\nc{\pp}{\perp}
\nc{\order}{{\cal O}}
\nc{\Eplus}{E_+}
\nc{\Eminus}{E_-}
\nc{\Epm}{E_\pm}
\nc{\Vav}{V_{\rm av}}
\nc{\Rin}{R_{\rm in}}
\nc{\Rout}{R_{\rm out}}
\nc{\eplus}{e_+}
\nc{\eminus}{e_-}
\nc{\epm}{e_\pm}
\nc{\eps}{\epsilon}
\nc{\veps}{\varepsilon}
\nc{\vnabla}{{\vec\nabla}}
\nc{\G}{\Gamma}
\nc{\w}{\omega}
\nc{\mh}{h}
\nc{\mg}{g}
\nc{\vphi}{\varphi}
\nc{\tlambda}{\tilde\lambda}
\nc{\be}{\begin{equation}}
\nc{\ee}{\end{equation}}
\nc{\ba}{\begin{eqnarray}}
\nc{\ea}{\end{eqnarray}}

\nc{\g}{\gamma}
\nc{\ol}{\overline}

\newtheorem{theo}{Theorem}[section]

\newtheorem{prop}{Proposition}[section]
\newtheorem{lem}{Lemma}[section]

\newtheorem{rmk}{Remark}[section]

\nc{\pT}{\partial_T}
\nc{\pz}{\partial_z}
\nc{\pt}{\partial_t}
\nc{\la}{\langle}
\nc{\ra}{\rangle}
\nc{\infint}{\int_{-\infty}^{\infty}}
\nc{\halfwidth}{6.5cm}
\nc{\figwidth}{10cm}

\nc{\nlayers}{L}
\nc{\nsectors}{M}
\nc{\indicator}{\mathbf{1}}
\nc{\Rhole}{R_{\rm hole}}
\nc{\Rring}{R_{\rm ring}}
\nc{\neff}{n_{\rm eff}}
\nc{\Frem}{F_{\rm rem}}
\nc{\Real}{\mathbb R}
\nc{\Z}{\mathbb Z}
\nc{\DD}{\Delta}
\nc{\cD}{\mathcal D}

\date{\today}

\begin{document}

\title{Periodic orbits in outer billiard}

\author{ {\em Alexander Tumanov}  and {\em Vadim Zharnitsky}\footnote{ 
Supported by NSF grant  DMS-0505216.} \\ \\
Department of Mathematics \\ 
University of Illinois at Urbana-Champaign \\ Urbana, IL 61801 \\
{\tt tumanov@math.uiuc.edu vz@math.uiuc.edu}
}

\maketitle

\begin{abstract}
It is shown that the set of 4-period orbits in outer billiard with piecewise
smooth convex boundary has an empty interior, provided that no four corners of
the boundary form a parallelogram.
\end{abstract}

\noindent

\section{Introduction}
The study of periodic orbits has been always important  in the field of 
Hamiltonian dynamics and  classical billiard is one of the early  
examples of a Hamiltonian dynamical system. This system was  introduced by  
G.D. Birkhoff, see {\em e.g.}  \cite{birk_1,birk_2}, who also showed that 
classical billiard with smooth convex boundary possesses at least 
two periodic orbits  of each $(p,q)$ type, see {\em e.g.} \cite{hassel}, 
for the proof of this result.

More recently, additional interest for the study of sets of periodic 
orbits came from the spectral theory of Laplace operator on bounded domains. 
V. Ivrii showed \cite{ivrii} that the so-called Weyl's asymptotics of 
distribution of large  eigenvalues of the Laplacian 
(with the Dirichlet or Neumann 
boundary conditions) 
holds,  provided periodic orbits of the corresponding classical 
system constitute the set of zero measure. Then, in order to fill the 
gap in the proof of what is sometimes called Weyl's conjecture \cite{weyl}, 
one has to demonstrate that the union of periodic orbits is 
a set of measure zero in the billiard phase space. Therefore, in contrast to   
Birkhoff theorem and its generalizations (see {\em e.g.} \cite{farber}), 
here one has to study the upper bound on the number 
(or rather measure) of periodic orbits.

Establishing that periodic orbits have zero measure turns out to 
be surprisingly 
hard for 
periodic orbits of arbitrary period. The only case which can be easily 
dealt with 
is the case of two period orbits. Indeed, 
since the segments of two period orbits must be normal to the boundary, then 
from  each  boundary point there can emanate only one 2-period orbit. 
These orbits then 
form at most one-parameter family which obviously has zero measure in the 
two dimensional phase space.

The case of three period orbits is already non-trivial. It was solved by 
Rychlik, see 
 \cite{rychlik}. The proof involved some symbolic calculations that were 
later removed 
in \cite{stojanov}. A much simpler proof relying on Jacobi's fields appeared
in \cite{wojtk}. Later in \cite{vorobetz} this result for three period orbits 
was extended to higher  dimensional billiards. All these proofs have been 
obtained by first demonstrating
that there are no open sets of periodic orbits and then verifying that 
the sets of positive measure do not exist either. The second part of the proof 
is relatively easy as the sets of positive measure have density points whose 
neighborhoods are ``almost'' open sets. 

More recently, another approach based on the theory of exterior differential 
systems  (EDS) has been proposed to study open sets of periodic orbits 
by Landsberg, 
Baryshnikov and the second author \cite{landsberg}, \cite{varna}. For systematic exposition 
of the EDS theory along with many applications, see {\em e.g.} \cite{ivey}.
Similar billiard formulation has been independently 
developed by T\" or\" ok \cite{torok}. The EDS approach gives 
a systematic  proof in the three period case  and reduces 4-period case to 
the study of zeros of certain system of polynomials. Unfortunately, the system 
is too hard to resolve even with the aid of symbolic calculations 
(at least by direct use of Maple or Mathematica).

In this paper, we consider a closely related system of outer (or dual) 
billiard, which 
is another popular model in Hamiltonian dynamics. Originally  introduced  
by Bernhard 
Neumann, the outer billiard was popularized by J\" urgen Moser \cite{moser}
and others as a model stability problem. 
See also the survey  article \cite{tabak} for more information and recent 
results on outer billiards.

The dynamics of the outer billiard is defined in the exterior of a convex boundary  
$\Gamma \in {\mathbb R}^2$ as follows: draw a line $L$ through a point 
$z_0=(x_0,y_0)\in {\mathbb R}^2$ tangent to $\Gamma$ in, say, the counterclockwise 
direction.  Find a point
$z_1=(x_1,y_1) \in L$ and such that the tangency point is dividing the segment
$|z_0,z_1|$ in half. The induced map $P: (x_0,y_0) \rightarrow (x_1,y_1)$
defines the outer (dual) billiard dynamics. The map is not well defined 
for lines for which 
 the tangency point is not unique. However such lines are countable and 
therefore,
the outer billiard map is not well defined on at most a set of zero measure 
\cite{dogru}. 
The exterior of the boundary can be then considered as a phase space and 
the invariant measure is given 
by the area form $\mu = dx \wedge dy$. 

The natural extension of the conjecture for classical billiard 
is that periodic orbits in outer billiard constitute the set of zero measure 
(except, may be, for some special boundaries, see section \ref{sec:special}). 
 While the outer billiard does not have such significance for the spectral 
asymptotics problem, our hope is that this study will help resolve 
related problem for the classical billiard. We have been also motivated by a 
recent article by Genin and Tabachnikov \cite{genin} 
which (among other results) contains a  proof that the set of 3-period orbits in outer billiard has an empty interior.

In this article, we study the set of  4-period orbits in outer billiards.
Our main result is contained in
\begin{theo}\label{theo:main}
Let $\Gamma \subset \Real^2$ be a piecewise smooth  convex closed curve. 
Assume that no four corners of $\Gamma$ form a parallelogram.
Then the set of 4-period orbits in the outer billiard has an empty interior.
\end{theo}
This theorem follows from Theorem \ref{theo:aux} on certain  
properties of the exterior differential system 
associated with the outer billiard. We recall that in the EDS approach, 
instead of 
asking which outer billiard boundaries possess open sets of 
periodic orbits, one studies which 2-parameter families of quadrilaterals 
can (or cannot) be orbits in an outer billiard. More precisely, 
we search for 2-dimensional  disks of quadrilaterals in the space of all 
quadrilaterals. These  2-dimensional disks must satisfy certain differential 
relations. 

In the next section,  we ``translate'' the problem in the language of exterior 
differential systems. This EDS corresponds to the Birkhoff distribution 
in the case of classical billiard \cite{landsberg,bz} and to the dual Birkhoff 
distribution in case of outer billiard \cite{tabak}. 
Then we find the solutions  of that EDS under some nondegeneracy conditions. 
As it turns out, for each 
nondegenerate quadrilateral there exists only one EDS solution, which 
corresponds to a 2-parameter family of 4-period orbits in an outer billiard. 
We verify that there are no other solutions by 
proving that Frobenius type integrability condition does not hold.

\section{EDS associated with outer billiard}
Since the set of initial conditions for which the billiard map is not well 
defined has zero measure,
we restrict ourselves to the complementary subset where the map is well 
defined. Below, we always assume that the outer billiard map is well defined.

We start with the proposition which establishes relation between open sets 
of $n$-periodic orbits and integral submanifolds in an associated exterior 
differential systems (see \cite{landsberg}, \cite{bz} or \cite{varna} for a related statement 
for the classical billiard).
\begin{prop}\label{prop:vz1}
Suppose that there exists an open set $Q$ of $n$-periodic orbits in 
the outer billiard phase space for 
the billiard with a  
convex piecewise smooth boundary $\Gamma \subset \Real^2$. Then there exists 
a 2-dimensional disk in the space of $n-$gons 
$M^2 \subset {\Real^{2 n}}$ such that 
\[
\theta^i|_{TM^2} = 0,
\]
where $i \in {\mathbb Z}/n{\mathbb Z}$ and 
\[
\theta^i = \frac{1}{2}(y_i-y_{i+1})d(x^{i}+x^{i+1}) -
\frac{1}{2} (x_i-x_{i+1})d(y^i+y^{i+1}).
\]
The following nondegeneracy conditions hold:
all points are different 
\[
(x_i,y_i)\neq (x_j,y_j) \,\,\, {\rm if} \,\,\, i \neq j 
\]
and no three consecutive points belong to the same line
\[
(x_{i-1}-x_i)(y_i -y_{i+1})\neq (x_i-x_{i+1})(y_{i-1}-y_i).
\]
The area form $dx^i\wedge dy^i \neq 0$ on $M^2$ for all $i$.
\end{prop}
{\em Proof:}
Let $M^2$ be the set of periodic orbits in ${\mathbb R}^{2n}$, the space 
of $n-$gons.
Any one-parameter family of periodic orbits $z_i(t)=(x_i(t),y_i(t))\in M^2$ , 
where $z_i\in \Real^2$ and $i\in {\mathbb Z}/n{\mathbb Z}$,  
must satisfy the tangency condition (middle point of any segment 
cannot move in the normal 
direction to the segment)
\[
\frac{d}{dt} \left ( \frac{z_{i+1}+z_i}{2}\right ) =\lambda (z_{i+1}-z_i),
\]
where $\lambda \in {\mathbb R}$. This relation 
implies that $\theta^i$ must vanish. 

In a sufficiently small neighborhood of each tangency point, the boundary is either smooth and convex or
it has a corner, therefore  
$M^2$ is an embedding of 
$Q$ in $\Real^{2n}$, the space of $n$-gons.
Since $Q$ is an open set, then $dx^1\wedge dy^1 \neq 0$.
Verification of nondegeneracy conditions is straightforward.
\hfill $\Box$ \\

By this proposition, it remains to find all two dimensional
integral submanifolds in the exterior differential system generated by
$\theta^i$ and satisfying the above nondegeneracy conditions.
We will refer to such integral manifolds as nondegenerate.

Next theorem gives the local description of the two dimensional integral 
manifolds in the outer billiard EDS.
\begin{theo}\label{theo:aux}
For  any nondegenerate convex 
quadrilateral 
there exists a unique nondegenerate 
connected integral 
manifold  containing the quadrilateral. This manifold is given by 
the quadrilaterals whose middle points coincide with those of the original 
quadrilateral. 
\end{theo}

\section{Proof of the theorem \ref{theo:aux}}

\subsection{New coframe}

Supplementing $\theta^i$ with $\omega^i$
\ba
\omega^i =  \frac{1}{2}(y_i-y_{i+1})\,d\,(x^{i}-x^{i+1}) -
\frac{1}{2} (x_i-x_{i+1})\, d\, (y^i-y^{i+1}), \label{eq:def_om}
\ea
we obtain a coframe $\{\theta^i,\omega^i  \}_{i=1}^{i=n}$.
It is easy to check that these $2n$ forms are linearly independent 
on an open dense subset of $\Real^{2n}:$ $\{ (x_i,y_i) \neq (x_j,y_j), 
(x_{i-1}-x_i)(y_i -y_{i+1})\neq 
(x_i-x_{i+1})(y_{i-1}-y_i),i,j = 1,2, ...,n \}$
using the following identities which can be directly verified
\ba
(y_i-y_{i+1})dx^{i+1} &-& (x_i-x_{i+1})dy^{i+1} = \theta^i-\omega^i \label{eq:id1}\\
(y_{i+1}-y_{i+2})dx^{i+1} &-& (x_{i+1}-x_{i+2})dy^{i+1} = \theta^{i+1}+\omega^{i+1}. \label{eq:id2}
\ea
Note that the determinant of the above linear system 
\[ \Delta_{i,i+1}=\left | \begin{array}{cc}
y_{i}-y_{i+1}  & - (x_{i}-x_{i+1})          \\
y_{i+1}-y_{i+2}  & - (x_{i+1}-x_{i+2})
\end{array} \right |   = 
\left | \begin{array}{cc}
x_{i}-x_{i+1}  & y_{i}-y_{i+1}          \\
x_{i+1}-x_{i+2}  & y_{i+1}-y_{i+2}
\end{array} \right |
\]
does not vanish by the nondegeneracy conditions 
(consecutive points do not belong to the same line).

This determinant $\DD_{i,i+1}$ has a clear geometric meaning.
It is the double area of the triangle with the vertices $(x_i,y_i), (x_{i+1},y_{i+1})$,
and $(x_{i+2},y_{i+2})$ (assuming the vertices are enumerated counterclockwise). 
The total area of the $n$-gon is an integral for the
system \cite{genin}. Indeed, adding the forms, we obtain
\ba
\sum_{i=1}^n \theta^i = \frac{1}{2} d \left ( \sum_{i=1}^n (y_i x_{i+1}-y_{i+1} x_i) \right ),
\ea
where the sum on the right handside is the total area of the n-gon.
 Therefore, for the quadrilateral $(n=4)$ 
we have
\ba
\DD_{1,2}+\DD_{3,4} = 2S ={\rm constant}\\
\DD_{2,3}+\DD_{4,1} = 2S ={\rm constant},
\ea
where $S$ is the area of the quadrilateral.

Solving (\ref{eq:id1}--\ref{eq:id2}), we obtain
\ba
dx^{i+1}=\frac{1}{\Delta_{i,i+1}} \left ( (x_{i+1}-x_{i+2})\omega^i + (x_{i}-x_{i+1})\omega^{i+1}
\right ) \label{eq:dxip1} \\
dy^{i+1}=\frac{1}{\Delta_{i,i+1}}\left ( (y_{i+1}-y_{i+2})\omega^i + (y_{i}-y_{i+1})\omega^{i+1}
\right ). \label{eq:dyip1}
\ea
In the last system and below  all relations are modulo 
the differential ideal generated by $\theta^i$.

\subsection{Exterior derivatives of the new coframe}

On the hypothetical integral manifold $M^2$, the differentials $d\theta^i$
must also vanish. Direct calculations show
\ba
d\theta^i =   dx^{i+1}\wedge dy^{i+1} - dx^i \wedge dy^i =0, i \in \Z/4\Z.
\label{eq:dtheta}
\ea
These identities are related to the area-conservation property of the outer billiard map. 

Another calculation gives the relation between some exterior products 
of the basis elements in the old and new coframes (by taking the exterior product of (\ref{eq:dxip1}) and 
(\ref{eq:dyip1}))
\ba
dx^{i+1} \wedge dy^{i+1} = -\frac{1}{\Delta_{i,i+1}}\omega^{i}\wedge
\omega^{i+1}.
\label{eq:area}
\ea
From (\ref{eq:dtheta}) and (\ref{eq:area}) we obtain that on $M^2$
the following relations hold
\ba
\Delta_{i,i+1}^{-1} \, \omega^{i}\wedge \omega^{i+1} =
\Delta_{i-1,i}^{-1} \, \omega^{i-1}\wedge \omega^{i}
\label{eq:rel}
\ea
for all $i$.

Now, we compute differentials of $\omega^i$:
\[
d \omega^i =
\frac{4}{\Delta_i}\, \omega^{i} \wedge \omega^{i+1},
\]
where we use the notation $\Delta_i := \Delta_{i,i+1}$.
Using (\ref{eq:rel}) we conclude that
\ba
d\omega^i = \frac{4}{\Delta_j}\, \omega^{j} \wedge \omega^{j+1}.
\label{eq:dom_eq}
\ea
for any $i,j \in \Z/4\Z$.

\subsection{The case of 3-period orbits}
Here we reproduce a result from \cite{genin} using EDS.
In this case $\DD_{i,i+1}$ is the double area enclosed by the triangular periodic 
orbit. Then the above relations simplify
\ba
\omega^i\wedge \omega^{i+1}=\omega^{i-1}\wedge \omega^i \\
d\omega^i = \frac{3}{S} \, \omega^j \wedge \omega^{j+1}, 
\ea
where $i\in \Z/3\Z ,j\in \Z/3\Z$ are arbitrary.

Since $\omega^1 \wedge \omega^2 \neq 0$ on $M^2$, then we must have 
a relation
\[
\omega^3 = a \, \omega^1 + b \, \omega^2.
\] 
Taking exterior product and using the above relations, we obtain
\[
\omega^1+\omega^2+\omega^3=0
\]
and therefore,
\[
d\omega^1+d\omega^2+d\omega^3=0 \Rightarrow 3\, \omega^1\wedge \omega^2 =0
\]
contradicting the independence of $\omega^1,\omega^2$.

\subsection{Integral elements}

On $M^2$ at most two 1-forms can be linearly independent.
Let us assume that $\omega^1 \wedge \omega^3 \neq 0$.
The case when $\omega^1\wedge \omega^3 = 0$ will be considered
separately.
The remaining 1-forms are then linearly dependent on $\omega^1,\omega^3$:
\ba
\omega^2 = a_1 \, \omega^1 + a_3 \, \omega^3 \\
\omega^4 = b_1 \, \omega^1 + b_3 \, \omega^3.
\ea
Taking the exterior product of both equations with $\omega^2,\omega^4$ and assuming 
${\mathcal D} = \Delta_2 \Delta_4 - \Delta_1 \Delta_3 \neq 0$
(the case ${\mathcal D}=0$ will be also evaluated separately), we obtain
\ba
0=\omega^2 \wedge \omega^2 &=& a_1\,  \omega^1 \wedge \omega^2 + a_3 \, \omega^3 \wedge \omega^2 \nn \\
0=\omega^4 \wedge \omega^4 &=& b_1\,  \omega^1 \wedge \omega^4 + b_3\,  \omega^3\wedge \omega^4 \nn \\
\omega^2 \wedge \omega^1 &=& a_3\,  \omega^3 \wedge \omega^1 \nn \\
\omega^4 \wedge \omega^1 &=& b_3 \, \omega^3 \wedge \omega^1.\nn
\ea
Using the relations in (\ref{eq:rel}), we obtain from the first two equations
\ba
\frac{a_1}{a_3}=\frac{\DD_2}{\DD_1} \,\,\, {\rm and} \,\,\, 
\frac{b_3}{b_1}=\frac{\DD_4}{\DD_3}  
\ea
and from the last two
\ba
\frac{a_3}{\DD_1}+\frac{b_3}{\DD_4}=0.
\ea
Expressing all the coefficients in terms of $b_3$ and then using the notation $v:= -b_3/\DD_4$, we 
find the relations
\ba
\omega^2 = \,\,\,\,\, v(\Delta_2\omega^1+\Delta_1 \omega^3) \label{eq:int_el1}\\
\omega^4 =-v (\Delta_3 \omega^1 + \Delta_4 \omega^3),
\label{eq:int_el2}
\ea
where  $v$ is a function defined on $M^2$.

Taking the exterior product of the above two equations, we obtain the relation
\ba
\omega^2\wedge \omega^4 = -v^2 {\mathcal D}\, \omega^1 \wedge \omega^3.
\label{eq:2413}
\ea

To compute $dv$, we need first to evaluate $d\Delta^i$. 
Using the 
definition of $\Delta_i$, the relations between the new and old coframes, 
and applying the following property 
\[\sum_{j=1}^4 \Delta_{i,j}= \sum_{j=1}^4 \left | \begin{array}{cc}
x_{i}-x_{i+1}  & y_{i}-y_{i+1}          \\
x_{j}-x_{j+1}  & y_{j}-y_{j+1}
\end{array} \right |   = \sum_{j=1}^4
\left | \begin{array}{cc}
x_{i}-x_{i+1}  & y_{i}-y_{i+1}          \\
0  & 0
\end{array} \right | = 0,
\]
we obtain 
\[
d \Delta^i = 
\frac{\DD_{i+2}}{\DD_{i-1}}\omega^i -\frac{\DD_{i}}{\DD_{i+1}}\omega^{i+2}+
\frac{\DD_i}{\DD_{i-1}}\omega^{i-1} -\frac{\DD_{i+2}}{\DD_{i+1}}\omega^{i+1}.
\]

Using the equations (\ref{eq:int_el1}--\ref{eq:int_el2}), we obtain
\ba
d \DD^1 = 
\frac{\DD_3}{\DD_4}(1-v(\DD_1+\DD_4))\omega^1 +
\frac{\DD_1}{\DD_2}(-1-v(\DD_2+\DD_3))\omega^3
\label{eq:ddelta1}\\
d \DD_2 = 
\frac{\DD_2}{\DD_1}(1+v(\DD_1+\DD_4))\omega^1 +
\frac{\DD_4}{\DD_3}(-1+v(\DD_2+\DD_3))\omega^3.
\label{eq:ddelta2}
\ea
Solving for $\omega^1,\omega^3$ in (\ref{eq:int_el1}-\ref{eq:int_el2}),
\ba
\omega^1 &=& \frac{1}{v {\mathcal D}} (\Delta_4 \omega^2 + \Delta_1 \omega^4) \\
\omega^3 &=& -\frac{1}{v {\mathcal D}} (\Delta_3 \omega^2 + \Delta_2 \omega^4)
\ea
and substituting these expressions in (\ref{eq:ddelta1}-\ref{eq:ddelta2}), 
we obtain
\ba
d\Delta^1 &=& \frac{\Delta_3}{\Delta_2 {\mathcal D} v}
(\Delta_1+\Delta_2 -{\mathcal D}v)  \, \omega^2 +
\frac{\Delta_1}{\Delta_4 {\mathcal D} v}(\Delta_3+\Delta_4+ {\mathcal D}v) \,  \omega^4  
\label{eq:del_om1}\\
d\Delta^2 &=& \frac{\Delta_4}{\Delta_1 {\mathcal D} v}
(\Delta_1+\Delta_2 +{\mathcal D}v)  \, \omega^2 +
\frac{\Delta_2}{\Delta_3 {\mathcal D} v}(\Delta_3+\Delta_4 -  {\mathcal D}v)  \, \omega^4. 
\label{eq:del_om2}
\ea
We can also express $\omega^2,\omega^4$ through $d\Delta^1,d\Delta^2$
by inverting the above equations:
\ba
\omega^2 &=-&\frac{1}{8S}
\left [ 
\frac{\Delta_2}{\Delta_3}(\Delta_3 + \Delta_4 -{\mathcal D}v) \, d\Delta^1 -
\frac{\Delta_1}{\Delta_4}(\Delta_3 + \Delta_4 +{\mathcal D}v) \, d\Delta^2
\right ] \label{eq:om_delta1} \\
\omega^4 &=& \frac{1}{8S}
\left [ 
\frac{\Delta_4}{\Delta_1}(\Delta_1 + \Delta_2 +{\mathcal D}v) \,  d\Delta^1 -
\frac{\Delta_3}{\Delta_2}(\Delta_1 + \Delta_2 -{\mathcal D}v) \,  d\Delta^2
\right ]. 
\label{eq:om_delta2}
\ea

\subsection{Special solutions}
\label{sec:special}
There exist  outer billiards in which  4-period orbits constitute a set of positive measure 
\cite{tab_private}. 
More precisely, the following statement holds:
\begin{prop}\label{prop:square}
Let $z_1,z_2,z_3,z_4$ be a convex quadrilateral and let 
$\zeta_i = (z_i +z_{i+1})/2, i= 1,2,3,4$ be an outer billiard defined by
the middle points of the initial quadrilateral. Then, there exist an open 
neighborhood $O(z_1) \subset \Real^2$, containing only 4-period points.
\end{prop}
{\em Proof:}
By the well known property of the triangle: midsegment between two sides 
is parallel to the 
third side, we observe that midpoints $\zeta_1, ..., \zeta_4$ form 
a parallelogram. 
Then, $z_1,..,z_4$ is a 4-period orbit in the 
outer billiard with this  boundary. It remains to show that if $z$ is 
sufficiently close to $z_1$ then $z$ is a footpoint of a 4-period orbit.
It follows easily from the same midsegment theorem assuming that we take 
small enough neighborhood so that its images do not intersect any of the lines 
containing the parallelogram sides.

\hfill $\Box$ \\

Therefore, each nondegenerate quadrilateral belongs to a two dimensional 
integral submanifold (actually it is a linear integral subspace) in 
the exterior differential system. Consider, the following specific example:
let the outer billiard be given by a unit square with vertices in 
$(0,0), (0,1), (1,0), (1,1)$. Let $z_1-(1/2,-i/2)$ be small. Then, the other 
vertices of periodic orbits are given by
\ba
z_1+z_2=2, z_2+z_3 = 2+ i 2, z_3+z_4= i2, z_1+z_4 = 0.
\ea
Solving this linear system, we obtain
\ba
x_4 & = & -x_1 \,\,\,\, x_2= 2-x_1 \,\,\,\, x_3  =  x_1 \nn \\ 
y_4 & = & -y_1 \,\,\,\, y_2=-y_1   \,\,\,\,\,\,\,\,\,\,\,\, y_3  =  y_1+2. \nn
\ea
Using (\ref{eq:def_om}) and the definition of $\Delta_{i}$, it is easy to compute:
\begin{eqnarray*}
\Delta_1 &=&  4(1-x_1) \\
\Delta_2 &=& 4(y_1+1) \\
\omega^1 &=& 2y_1 dx^1 - 2(x_1-1)dy^1 \\
\omega^2 &=& 2(y_1+1)dx^1 + 2(1-x_1)dy^1 \\
\omega^3 &=& 2(y_1+1)dx^1 - 2x_1 dy^1
\end{eqnarray*}
and then substituting these expressions in (\ref{eq:int_el1}),
we obtain 
\ba
v=\frac{1}{4(y_1+1-x_1)} \Rightarrow v=\frac{1}{\Delta_2-\Delta_3}.
\label{eq:spec_rel}
\ea

This calculation shows that for $v=1/(\Delta_2-\Delta_3)$ there exists 
a solution for each quadrilateral\footnote{Indeed, since outer billiard 
map commutes with affine transformations, similar solution with 
the same relation (\ref{eq:spec_rel}) exists for arbitrary nondegenerate 
convex quadrilateral. } and  there are no other solutions
with such $v$. Indeed,  this modified EDS has 2 additional 1-forms which 
must vanish and which are linearly independent:
\ba
\theta^5  & = & \DD_2 \, \omega^1 + (\DD_3-\DD_2)\, \omega^2 + \DD_1\, 
\omega^3 \\
\theta^6  & = & \DD_3 \, \omega^1 + \DD_4 \, \omega^3 + (\DD_2-\DD_3)\, 
\omega^4.\, 
\ea
But then, by the standard ODE argument there is at most one solution.

\subsection{Computation of $du$}
Knowing the special solution, it is now convenient to change the parameter
\[
v = \frac{u}{\DD_2-\DD_3},
\]
so that the special solution corresponds to $u\equiv 1$.
Taking exterior derivative of (\ref{eq:int_el1}) with $v$ replaced by 
$u/(\DD_2-\DD_3)$
\ba
(\DD_2-\DD_3)\omega^2 = u(\DD_2 \omega^1 + \DD_1 \omega^3),
\ea 
we obtain (using $d\DD_1=-d\DD_3$ and $d\omega^i=d\omega^j$)
\ba
d(\DD_1+\DD_2) \omega^2 + (\DD_2-\DD_3) d \omega^2 = \nn \\
\frac{\DD_2-\DD_3}{u} du \wedge \omega^2 +
u(d\DD^2 \wedge \omega^1 + d\DD^1 \wedge \omega^3)+
u(\DD_2  + \DD_1) d\omega^2.
\label{eq:domega2}
\ea
Now, using (\ref{eq:del_om1}-\ref{eq:del_om2}) and (\ref{eq:dom_eq}), we obtain
\ba
\left [ \frac{\DD_1}{\DD_4 {\cD} v}(\DD_3 +\DD_4+\cD v) + 
\frac{\DD_2}{\DD_3 \cD v }(\DD_3+\DD_4-\cD v)   \right ] \w^4\wedge w^2 +
d\w^2 (\DD_2-\DD_3-u(\DD_1+\DD_2)) =  \nn \\ 
\frac{\DD_2-\DD_3}{u} du \wedge \omega^2 +
u \left [ \frac{\DD_4}{\DD_3}(-1 + v(\DD_2+\DD_3))-
\frac{\DD_3}{\DD_4}(-1 + v(\DD_1+\DD_4)) \right ]\w^3\wedge \w^1. \nn
\ea
The last expression can be further simplified
\ba
\cD du\wedge \w^2 + 4(\DD_1+\DD_2-S)(1-u)\, \w^2\wedge \w^4 = 0,
\label{eq:du2}
\ea
where we used:
\[
\cD=\DD_2 \DD_4 - \DD_1 \DD_3 = (\DD_1-\DD_2)(\DD_2-\DD_3) 
\]
\[
\cD v = (\DD_1-\DD_2)u 
\]
\[
d\w^i=\frac{4}{\DD_2}\w^2\wedge \w^3 = 4 v \w^1 \wedge \w^3 = 
-\frac{4}{v \cD}\w^2\wedge \w^4 = -\frac{4}{u(\DD_1-\DD_2)} \w^2\wedge \w^4 
\]
\[
\DD_1+\DD_3 = \DD_2+\DD_4 = 2S \,\, (={\rm constant}).
\]
To derive similar relation for $du\wedge \w^4$, we add up
(\ref{eq:int_el1}-\ref{eq:int_el2})
\[
(\DD_2-\DD_3)(\w^2+\w^4)=u((\DD_2-\DD_3)\, \w^1+(\DD_1-\DD_4)\, \w^3)
\]
and since $\DD_1-\DD_4=\DD_2-\DD_3$, we have
\[
(\DD_2-\DD_3)\, (\w^2+\w^4)=u\, (\DD_2-\DD_3)\, (\w^1+\w^3) \Rightarrow
\w^2+\w^4 = u\, (\w^1+\w^3).
\]
Taking exterior derivative, we obtain
\[
2d\\, w^i = du\wedge (\w^1+\w^3) + 2u \, d\w^i,
\]
which implies
\[
-\frac{8}{u(\DD_1-\DD_2)}\, \w^2\wedge\w^4 = \frac{1}{u} du\wedge (\w^2+\w^4)-
\frac{8}{\DD_1-\DD_2} \, \w^2\wedge\w^4 
\]
and after multiplying with $u\cD$
\[
\cD du \wedge (\w^2+\w^4) +8(1-u)(\DD_2-\DD_3)\w^2\wedge\w^4.
\]
Subtracting (\ref{eq:du2}) from the last expression, we obtain
\ba
\cD \, du\wedge \w^4 + 4(\DD_2-\DD_3-S) (1-u) \, \w^2\wedge\w^4=0
\ea
with  (\ref{eq:du2}) rewritten in a similar form:
\ba
\cD \, du\wedge \w^2 + 4(\DD_2-\DD_3+S) (1-u) \, \w^2\wedge\w^4=0.
\ea
Then 
\ba
\cD\, du = (1-u)\left ( (S-\DD_2+\DD_3)\, \w^2 + (\DD_2-\DD_3 +S)\, \w^4 \right ) .
\ea

Substituting (\ref{eq:om_delta1}-\ref{eq:om_delta2}) in the last expression 
we obtain
\ba
\frac{8S}{1-u} \, du = (a_1 u +b_1)\, d \DD^1 + (a_2 u + b_2)\, d\DD^2, 
\label{eq:dud1d2}
\ea
where
\ba
a_1 &=&\frac{\DD_4}{\DD_1} \left ( +1+\frac{S}{\DD_2-\DD_3} \right )
+ \frac{\DD_2}{\DD_3} \left ( -1+\frac{S}{\DD_2-\DD_3} \right ) \\
a_2&=&\frac{\DD_1}{\DD_4} \left ( -1+\frac{S}{\DD_2-\DD_3} \right )
+ \frac{\DD_3}{\DD_2} \left ( +1+\frac{S}{\DD_2-\DD_3} \right ) \\
b_1&=&\frac{1}{\cD}\left [ \frac{\DD_4}{\DD_1}(\DD_1+\DD_2)(S+\DD_2-\DD_3) 
-\frac{\DD_2}{\DD_3}(\DD_3+\DD_4)(S-\DD_2+\DD_3) \right ] \\
b_2&=&\frac{1}{\cD}\left [ \frac{\DD_1}{\DD_4}(\DD_3+\DD_4)(S-\DD_2+\DD_3) 
-\frac{\DD_3}{\DD_2}(\DD_1+\DD_2)(S+\DD_2-\DD_3) \right ]. 
\ea
Taking the exterior derivative of (\ref{eq:dud1d2}) we obtain
\ba
0= d\left ( \frac{8S}{1-u} du \right ) = \left  (
-(u \D_2 a_1 + \D_2 b_1 + a_1 \D_2 u)\,  +
(u \D_1 a_2 + \D_1 b_2 + a_2 \D_1 u) \right ) d\DD^1 \wedge d\DD^2, 
\ea
where $\D_i := \frac{\D}{\D \DD_i}$. 

Next, we use the expression for $du$, to replace $\D_i u$ 
by $(1-u)(a_i u+ b_i)/8S$:
\ba
u \left ( \D_2 a_1 -\D_1 a_2 + \frac{a_2 b_1 - a_1 b_2}{8 S} \right )
+ \D_2 b_1 - \D_1 b_2 + \frac{a_2 b_1 - a_1 b_2}{8 S} = 0.
\ea
Using Maple and then some simplifications, we compute
\begin{eqnarray*}
\frac{S(u-1)}{\DD_1 \DD_2 \DD_3 \DD_4 \cD}\,\,
( - \DD_1^{4} + 4\, \DD_1^3\,S + 5\, \DD_1^2\, \DD_2^2 - 10\,\DD_1^2\,
\DD_2\,S - 3\,\DD_1^2\,S^2 + 20\,\DD_1\, \DD_2\,S^2 -
\end{eqnarray*}
\begin{eqnarray*}
  2\, \DD_1\,S^3 - 10\, \DD_1\, \DD_2^2\,S - 2\,S^3\, \DD_2 - 
 \DD_2^4 + 4\, \DD_2^3\,S - 3\, \DD_2^2\,S^2) = 0.
\end{eqnarray*}
The last equality cannot hold identically on an open subset. Indeed,  
we have imposed $u\neq 1$ and the numerator is a nontrivial polynomial
in two variables and thus, cannot vanish on an open set. 

\begin{rmk}\label{rmk:last}
It is interesting to note that the last expression has the form
$(u-1)f(\DD_1,\DD_2)$. In other words, $u-1$ can be factored out 
again making
the calculations much easier. In general, one might expect 
this expression to be of the form: $f(\DD_1,\DD_2) u + g(\DD_1,\DD_2)=0$.
Then we would have to take the exterior derivative once more and then check 
solvability condition  of this new system.
\end{rmk}

\subsection{Degenerate cases: $\omega^1\wedge \omega^3 =0$ or ${\mathcal D}=0$.}\subsubsection{$\omega^1\wedge\omega^3=0$}
\begin{lem}\label{lem:d1d4}
\ba
\omega^1\wedge \omega^3 = 0 \,\,\,\, {\rm implies} \,\,\,\,  \Delta_1=\Delta_2
\,\,\, {\rm or} \,\,\, \Delta_1=\Delta_4.
\label{eq:o1o3}
\ea
\end{lem}
{\em Proof:} \\
Since $\omega^1\wedge \omega^2 \neq 0$ we can represent integral elements by
\ba
\omega^3 = a_1 \omega^1 + a_2 \omega^2 \label{eq:int_deg1}\\
\omega^4 = b_1 \omega^1 + b_2 \omega^2
\label{eq:int_deg2}
\ea
and then
\[
\omega^3 \wedge \omega^1 = a_2 \omega^2 \wedge \omega^1 \Rightarrow a_2 = 0
\]
and
\[
\omega^3\wedge \omega^2 = a_1 \omega^1 \wedge \omega^2.
\]
Using the relations (\ref{eq:rel}), we then obtain
\[
a_1 = -\frac{\Delta_2}{\Delta_1}.
\]
Therefore, we have
\[
\Delta_1 \omega^3 + \Delta_2 \omega^1  = 0.
\]
Similarly, we obtain for $b_2$,
\[
\omega^4 \wedge \omega^1  = b_2 \omega^2 \wedge \omega^1,
\]
which with (\ref{eq:rel}) implies
\[
b_2 = -\frac{\Delta_4}{\Delta_1}.
\]
For $b_2$, we also have
\[
\omega^4 \wedge \omega^3 = b_2 \omega^2 \wedge \omega^3.
\]
Then, using (\ref{eq:rel}) once more, we have
\[
b_2 = -\frac{\Delta_3}{\Delta_2}.
\]
Now, using both equations for $b_2$, we obtain
\[
\Delta_2 \Delta_4 = \Delta_1 \Delta_3 \Rightarrow {\mathcal D}=0.
\]

On the other hand, ${\mathcal D} = 0$ implies 
\ba
\cD = (\Delta_1-\Delta_2)(\Delta_1-\Delta_4)=0.
\label{eq:Dzero}
\ea
Therefore, in some neighborhood of $M^2$, either $\Delta_1 = \Delta_2$ or
$\Delta_1=\Delta_4$.

{\hfill $\Box$}

Suppose, first that $\omega^2\wedge \omega^4 = 0$,
then
\[
\omega^1\wedge \omega^3 = \omega^2\wedge \omega^4 = 0.
\]
Taking exterior products of (\ref{eq:int_deg1}-\ref{eq:int_deg2}) with
$\omega^i$ and using (\ref{eq:rel}), we obtain
\ba
\omega^3 = - \frac{\Delta_2}{\Delta_1}\omega^1 \\
\omega^4 = -\frac{\Delta_4}{\Delta_1}\omega^2.
\ea
However, either $\Delta_1 = \Delta_2$ or $\Delta_1=\Delta_4$. In the first
case, we obtain
\[
\omega^3+\omega^1 = 0 \Rightarrow d \omega^3= - d\omega^1,
\]
but this contradicts (\ref{eq:dom_eq}). Similarly, in the second
case ($\Delta_1=\Delta_4$), we obtain $d\omega^4 = -d\omega^2$
also leading to contradiction.

Now, we are left to consider the case $\omega^2\wedge \omega^4
\neq 0$. By relabeling, this case is equivalent to the case
considered in the next section ${\mathcal D}=0, \omega^2\wedge
\omega^4 = 0, \omega^1\wedge \omega^3 \neq 0$. 

\subsection{${\mathcal D}=0$, $\omega^1\wedge\omega^3 \neq 0$}
In this case we can use representation of integral elements
(\ref{eq:int_el1}-\ref{eq:int_el2}) and then we also have
$\omega^2\wedge \omega^4=0$.
Using (\ref{eq:Dzero}) we also have that $\Delta_1=\Delta_2$ or
$\Delta_1=\Delta_4$. Assume first that
 \ba
 \Delta_1=\Delta_2 \Rightarrow \Delta_3 = \Delta_4.
\ea

Now, using formulae (\ref{eq:ddelta1}-\ref{eq:ddelta2}), we obtain
\[
v(\Delta_1+\Delta_4)=0,
\]
which can only occur if $v=0$ on $M^2$. Then,
by (\ref{eq:int_el1}--\ref{eq:int_el2}) two 1-forms vanish
identically $\omega^2=\omega^4 = 0$ implying $\omega^1\wedge\omega^2 = 0$,
which contradicts the genericity assumption.

In the second scenario $\Delta_1 = \Delta_4$, similar
calculations lead to the same contradiction. \\ \\

\noindent
\begin{large}
{\bf Acknowledgment}
\end{large}
We would like to thank S. Tabachnikov for pointing out that there are 
outer billiards with two parameter families of 4-period orbits.

\end{document}